\documentclass[11pt,reqno]{amsart}
\usepackage{amssymb,amsmath,tabularx}
\usepackage{amsthm,verbatim}
\usepackage[all]{xy}
\usepackage[bookmarks=true]{hyperref}
\usepackage{pstricks,pst-node,pst-plot}
\numberwithin{equation}{section}

\newtheorem{thm}{Theorem}[section]
\newtheorem{lem}[thm]{Lemma}
\newtheorem{prop}[thm]{Proposition}
\newtheorem{cor}[thm]{Corollary}
\theoremstyle{definition}

\newtheorem{rem}[thm]{Remark}
\newtheorem{ques}[thm]{Question}

\newcommand{\Res}{\operatorname{Res}}
\newcommand{\Ind}{\operatorname{Ind}}
\newcommand{\sgn}{\operatorname{sgn}}
\newcommand{\sym}[1]{\mathfrak{S}_{#1}}

\begin{document}
\title{Specht modules with abelian vertices}

\author{Kay Jin Lim}
\address{Department of Mathematics, National University of Singapore, Block S17, 10 Lower Kent Ridge Road, Singapore 119076.}
\email{matlkj@nus.edu.sg}

\date{January 2011}

\thanks{Supported by MOE Academic Research Fund R-146-000-135-112.}


\begin{abstract} In this article, we consider indecomposable Specht modules with abelian vertices. We show that the corresponding partitions are necessarily $p^2$-cores where $p$ is the characteristic of the underlying field. Furthermore, in the case of $p\geq 3$, or $p=2$ and $\mu$ is $2$-regular, we show that the complexity of the Specht module $S^\mu$ is precisely the $p$-weight of the partition $\mu$. In the latter case, we classify Specht modules with abelian vertices. For some applications of the above results, we extend a result of M. Wildon and compute the vertices of the Specht module $S^{(p^p)}$ for $p\geq 3$.
\end{abstract}

\maketitle

\section{Introduction} The representation theory of symmetric groups has been one of the major research areas since early on in the last century. The theory has been well-developed but yet it seems that very little is known for the modular case. One way of understanding the structure of representations of finite groups is through the notion of relative projectivity, on which J. A. Green \cite{JG1} defined the vertices of modules of finite groups about fifty years ago. The vertices of modules are, in some way, related to the complexity of the modules defined by J. Alperin and L. Evens \cite{JALE}, and the rank variety defined by J. Carlson \cite{JC}.

Classically, Young modules, Specht modules and simple modules are the major objects extensively studied in the representation theory of symmetric groups. The computation of the vertices of Young modules has been done by J. Grabmeier \cite{JGrab}. The vertices of signed Young modules have been computed by S. Donkin \cite{SD}. However, the vertices of simple modules and Specht modules remain mostly unknown. The vertices of Specht modules of hook shape in $p=2$ case were first considered by G. M. Murphy and M. H. Peel in \cite{GMMP}, with a mistake which has been corrected by M. Wildon \cite{MW}. M. Wildon has made some progress on the computation of the vertices of Specht modules which he has computed the vertices of simple Specht modules of hook shape \cite{MW} and showed that, in general, a vertex of a Specht module contains some large $p$-subgroup \cite{MW2}. Recently, there are some computation of the vertices of simple modules made by S. Danz, B. K\"{u}lshammer, J. M\"{u}ller and R. Zimmermann \cite{SDBKRZ,JMRZ}.

It is well known that, in general, if the defect groups of an indecomposable module are abelian then its vertices are necessarily abelian too. In particular, for the representations of symmetric groups in the modular case, representations with abelian defect correspond to partitions of $p$-weights strictly less than the characteristic $p$ of the underlying field. But there are examples of Specht modules of hook shape whose vertices are abelian but yet their defect groups are not abelian. For example, for $p=3$ the Specht module $S^{(7,1^3)}$ has vertices the Sylow $3$-subgroups of $\sym{6}\times \sym{3}$ (see \cite[Theorem 2]{MW}) but with defect groups the Sylow $3$-subgroups of $\sym{9}$.

In this article, we shall mainly be concerned with Specht modules with abelian vertices. Motivated by the question of classifying these Specht modules, we give some necessary conditions for such Specht modules. Indeed, for $p\geq 3$, we show that no other abelian subgroups can be the vertices of Specht modules besides the elementary abelian ones. For $p=2$, we could not give a definite answer as in the odd characteristic case unless the corresponding partitions are $2$-regular. Under the hypothesis that the vertices of the Specht modules are abelian, a class of partitions arises naturally; namely the $p^2$-core partitions. In fact, for $p=2$, we show that a $2$-regular partition $\mu$ is a $4$-core if and only if the Specht module $S^\mu$ has elementary abelian vertices. For all the Specht modules mentioned above, we conclude that their complexities are precisely the $p$-weights of their labeled partitions.

We organize the article in the following way. In \S \ref{S: notation}, we lay down some basic knowledge about the representations of symmetric groups, the complexities and the rank varieties for modules. We state our main results in \S \ref{S: main results} and prove them in \S \ref{S: proofs}. In \S \ref{S: consequences}, we draw some consequences of the main results, in which we generalize a result of M. Wildon and show that the vertices of $S^{(p^p)}$ are the Sylow $p$-subgroups of $\sym{p^2}$ when $p\geq 3$. This is an example where the partition $(p^p)$ is a $p^2$-core but the Specht module $S^{(p^p)}$ does not have abelian vertices. In the last section \S \ref{S: some further questions}, we post some questions which arise naturally from our results.

\section{Preliminaries}\label{S: notation}

We introduce the notation and background which we require. General references for this section are \cite{DB,GJAK,GJ1}.

\subsection{The representations}\label{S: the representations} Let $G$ be a finite group and $F$ be a field of characteristic $p$. In this article, all $FG$-modules are finite dimensional vector spaces over $F$.

Let $M$ be a $FG$-module and $H$ be a subgroup of $G$. We regard the restriction $\Res^G_H M$ as the $FH$-module in the obvious way. If $S$ is a $FH$-module we write $\Ind^G_H S$ for the induced module. Let $N$ be another $FG$-module. We write $N\mid M$ if $N$ is isomorphic to a direct summand of $M$ as $FG$-modules.

Suppose that $p>0$. We say that the module $M$ is relatively $H$-projective if $M\mid \Ind^G_H \Res^G_H M$. A vertex $Q$ of an indecomposable $FG$-module $M$ is a minimal subgroup of $G$ subject to the condition that $M$ is relatively $Q$-projective. Given a vertex $Q$ of $M$, a source of $M$ is an indecomposable $FQ$-module $S$ such that $M\mid \Ind^G_Q S$.

Let $FG=I_1\oplus\cdots\oplus I_m$ be a decomposition of the $F(G\times G)$-module $FG$, with the action given by $(g,h)x=gxh^{-1}$, into indecomposables. Suppose that $1=\sum_{i=1}^m e_i$ with $e_i\in I_i$ for each $1\leq i\leq m$. Indeed, $I_i=e_i FG$. The elements $e_i$ are mutually orthogonal primitive central idempotents of the algebra $FG$ and they are called the blocks of $FG$. For an indecomposable $FG$-module $M$, there is a unique block $e_j$ such that $e_jM\neq 0$ and $e_iM=0$ for all $i\neq j$. In this case, $e_jM=M$ and we say that the module $M$ lies in the block $e_j$. Note that $FG\cong \Ind^{G\times G}_{\Delta(G)} F$ where $\Delta(G)$ is the diagonal embedding of $G$ into $G\times G$. Thus, for each $1\leq i\leq m$, a vertex of $I_i$ is of the form $\Delta(D_i)$ for some subgroup $D_i$ of $G$ and the subgroup $D_i$ is called a defect group of the block $e_i$.

\begin{thm}[\cite{JG1}] Let $F$ be a field of characteristic $p>0$, $G$ be a finite group and $M$ be an indecomposable $FG$-module. Then
\begin{enumerate}
  \item [(i)] any vertex of $M$ is a $p$-subgroup of $G$,
  \item [(ii)] the vertices of $M$ are conjugate to each other in $G$, and
  \item [(iii)] if $M$ lies inside the block $e_j$ of $FG$ then $M$ is relatively $D_j$-projective where $D_j$ is a defect group of $e_j$. In particular, a vertex of $M$ is necessarily a subgroup of $D_j$ up to conjugation.
\end{enumerate}
\end{thm}

Suppose that $L$ is a field extension of $F$. We write $L\otimes_F M$ for the $LG$-module upon field extension. In the latter part of this article, we are required to deal with vertices of modules over field extension. We include a little lemma here.

\begin{lem}[{\cite[\S III Lemma 4.14]{WF}}]\label{L:vtx fld ext} Let $L$ be a field extension of $F$ and $H$ be a $p$-subgroup of a finite group $G$. Then an $FG$-module $M$ is relatively $H$-projective if and only if $L\otimes_F M$ is relatively $H$-projective.\\
In particular, if $M$ is indecomposable and $L\otimes_F M$ remains indecomposable as a $LG$-module then a $p$-subgroup $Q$ of $G$ is a vertex of $M$ if and only if $Q$ is a vertex of $L\otimes_F M$.
\end{lem}

\subsection{Rank varieties}\label{S:varieties} Let $p$ be a prime. The $p$-rank of a finite group $G$ is the largest integer $n$ subject to the condition that $G$ contains an elementary abelian $p$-subgroup of order $p^n$. The $p$-rank of the abelian $p$-group $\mathbb{Z}_{p^{n_1}}\times \cdots \times \mathbb{Z}_{p^{n_m}}$ with $n_1,\ldots,n_m>0$ is $m$.

Let $E$ be an elementary abelian $p$-group of $p$-rank $n$ with generators $g_1,\ldots,g_n$. For each non-zero point $\alpha=(\alpha_1,\ldots,\alpha_n)$ of the space $F^n$, we write $u_\alpha=1+\sum_{i=1}^n \alpha_i(g_i-1)\in FE$. Note that $\langle u_\alpha\rangle$ is a cyclic $p$-group. Suppose that $F$ is algebraically closed. The rank variety of a $FE$-module $M$ is the set $$V^\#_E(M)=\{\mathbf{0}\}\cup \{\mathbf{0}\neq \alpha\in F^n\,|\, \text{$M$ is not projective as $F\langle u_\alpha\rangle$-module}\}.$$ The rank variety $V^\#_E(M)$ is a homogeneous and closed subvariety of the affine space $\mathbb{A}^n(F)$ \cite[Theorem 4.3]{JC}. We shall write $\dim V^\#_E(M)$ for the dimension of the algebraic variety $V^\#_E(M)$. The rank variety depends on the choice and order of the generators of $E$, but its dimension does not. In the case where $p\nmid \dim_F M$, for each $\mathbf{0}\neq \alpha\in F^n$, $M$ necessarily has a summand of dimension coprime to $p$ as $F\langle u_\alpha\rangle$-modules. In this case, we have $V^\#_E(M)=V^\#_E(F)$ and thus $\dim V^\#_E(M)=n$. Let $N$ be another $FE$-module. It is clear from the definition that $V^\#_E(M\oplus N)=V^\#_E(M)\cup V^\#_E(N)$.

Let $G$ be a finite group and $M$ be a $FG$-module. Let $$\mathbf{P}:\cdots \to P_2\to P_1 \to P_0\to M$$ be a minimal projective resolution of $M$. Following \cite{JALE}, the complexity $c$ of the module $M$ is the smallest non-negative integer such that $$\lim_{n\to\infty} \frac{\dim_F P_n}{n^c}=0;$$ in other words, the polynomial rate growth of the dimensions of the terms of the sequence $\mathbf{P}$. We write $c_G(M)$ for the number $c$ to indicate the dependence of $c$ on $M$ as $FG$-module.

\begin{thm}\label{P:complexity of mod}\
\begin{enumerate}
\item [(i)] \textup{({\cite[\S 1 Theorem]{JALE}, \cite[Theorem 1.1]{GALS} and \cite[Theorem 7.6]{JC}})} Let $M$ be an $FG$-module, $\mathcal{E}$ be a set of representatives of elementary abelian $p$-subgroups of $G$ up to conjugation and $\mathcal{E}^{\text{max}}$ be a set of representatives of maximal elementary abelian $p$-subgroups of $G$ up to conjugation. Then $$c_G(M)=\max_{E\in\mathcal{E}}\{\dim V^\#_E(M)\}=\max_{E\in\mathcal{E}^{\text{max}}}\{\dim V^\#_E(M)\}.$$

\item [(ii)] \cite[Corollary 4,5]{JALE} The complexity of an indecomposable $FG$-module $M$ is bounded above by the $p$-rank of a defect group of the block in which $M$ lies.
\end{enumerate}
\end{thm}

Following Theorem \ref{P:complexity of mod} (i), we have the conclusion that $c_G(M)$ is the $p$-rank of $G$ provided $p\nmid \dim_F M$. We have a more refined statement of Theorem \ref{P:complexity of mod} (ii).

\begin{prop}[{\cite[Lemmas 5.2, 5.3]{JALE}}]\label{L:complx rel vtx} Let $M$ be a $FG$-module where $F$ is a field of characteristic $p>0$.
\begin{enumerate}
\item [(i)] Suppose that $M$ is relatively $H$-projective for some subgroup $H$ of $G$. Then $c_G(M)=c_H(M)$.
\item [(ii)] Suppose that $M$ is indecomposable. Let $Q$ be a vertex of $M$ with a source $S$. Then $c_G(M)=c_Q(S)$.
\end{enumerate}
\end{prop}

Readers who are familiar with the variety theory for modules would have realized that the author has avoided taking cohomological variety into discussion, thanks to the Quillen stratification theorem and \cite[Theorem 1.1]{GALS}.

\subsection{The representations of symmetric groups} We now briefly go through the representations of symmetric groups. Let $n$ be a non-negative integer. A partition $\mu$ of $n$ is a sequence of positive integers $(\mu_1,\mu_2,\cdots,\mu_s)$ such that $\mu_1\geq \mu_2\geq \cdots\geq \mu_s$ and $\sum_{i=1}^s\mu_i=n$. In this case, we write $|\mu|=n$. Note that we allow the empty partition $\varnothing$ to be the unique partition of $0$. The partition $\mu$ is called $p$-singular if there is some $i\geq 0$ such that $\mu_{i+1}=\mu_{i+2}=\cdots=\mu_{i+p}>0$; otherwise, it is called $p$-regular. Let $\Lambda(n)$ be the set consisting of all partitions of $n$. There is a one-to-one correspondence between $\Lambda(n)$ and the set of Young diagrams with $n$ nodes in an obvious way. The Young diagram of $\mu$ is written as $[\mu]$. Fix a positive integer $m$, not necessarily a prime. Each partition $\mu$ is associated to a partition $\widetilde{\mu}$ and an non-negative integer $w$ such that $|\mu|=|\widetilde{\mu}|+mw$. The partition $\widetilde{\mu}$ and the integer $w$ is called the $m$-core and the $m$-weight of $\mu$, respectively \cite[\S 2.7]{GJAK}. The $m$-core of $\mu$ is obtained by successively removing $w$ removable rim $m$-hooks. In the case where $w=0$ or equivalently $\mu=\widetilde{\mu}$, we say that the partition $\mu$ is an $m$-core.

We write $\sym{n}$ for the symmetric group acting on $n$ letters. To each partition $\mu$ of $n$, we have the $F\sym{n}$-module $S^\mu_F$, the Specht module labeled by the partition $\mu$. We usually write $S^\mu$ for $S^\mu_F$ if the underlying field is understood. The dimension of the Specht module $S^\mu_F$ is given by the hook formula $$\dim_F S^\mu_F=\frac{n!}{\prod_{(i,j)\in [\mu]} h_{i,j}}$$ where $h_{i,j}$ denotes the hook length of the node $(i,j)$ of $[\mu]$ and $n=|\mu|$. We note that the dimension of $S^\mu_F$ is independent of the field $F$.

\begin{thm}\label{P:Specht mod} Let $p$ be the characteristic of a field $F$ and $n$ be a positive integer.
\begin{enumerate}
  \item [(i)] \cite[2.7.40]{GJAK} Let $m$ be a positive integer. The number of hook lengths of a Young diagram $[\mu]$ which are divisible by $m$ is precisely the $m$-weight of $\mu$.
  \item [(ii)] \cite[Theorem 4.12]{GJ1} For $p=0$, the set $\{S^\mu\,|\,\mu\in\Lambda(n)\}$ is a complete set of non-isomorphic simple $F\sym{n}$-modules.
  \item [(iii)] \cite[Corollary 13.18]{GJ1} Let $\mu$ be a partition. For $p\geq 3$, or $p=2$ and $\mu$ is $2$-regular, the Specht module $S^\mu$ is indecomposable.
  \item [(iv)] \cite[\S 6.1 and Theorem 6.2.45]{GJAK} The blocks of $F\sym{n}$ are parametrized by the $p$-cores of the partitions of $n$ such that the block $e_{\widetilde{\mu}}$ labeled by the $p$-core $\widetilde{\mu}$ contains the Specht module $S^\mu$. In particular, if $\mu,\lambda$ are partitions of $n$ then the Specht modules $S^\mu, S^\lambda$ lie inside the same block of $F\sym{n}$ if and only if $\widetilde{\mu}=\widetilde{\lambda}$. A defect group of the block $e_{\widetilde{\mu}}$ is conjugate to a Sylow $p$-subgroup of $\sym{wp}$ where $w$ is the $p$-weight of $\mu$.
\end{enumerate}
\end{thm}

We shall take a step further to discuss the vertices of a special class of modules of the symmetric groups, the Young modules \cite{GJ2}. Here, and hereafter, whenever we have subgroups $H\subseteq \sym{a}$, $K\subseteq \sym{b}$ where $a+b\leq n$ for some positive integers $a,b,n$, we write $H\times K\hookrightarrow \sym{a}\times \sym{b}\hookrightarrow \sym{n}$ for the obvious inclusions.

Let $\mu=(\mu_1,\ldots,\mu_r)$ be a partition of $n$. Let $\mathfrak{S}_{\mu}$ be the Young subgroup of $\sym{n}$, i.e. $$\mathfrak{S}_\mu=\sym{\mu_1}\times \sym{\mu_2}\times\cdots\times \sym{\mu_r}.$$ Let $M^\mu\cong \Ind^{\sym{n}}_{\mathfrak{S}_\mu} F$ be the associated permutation module. Suppose that $M^\mu=M_1\oplus\cdots\oplus M_k$ is a decomposition of $M^\mu$ into indecomposable summands. Then there is a unique direct summand $M_j$ of $M^\mu$ such that $S^\mu\subseteq M_j$. We write $Y^\mu$ for $M_j$ and call it a Young module. Clearly, the Young modules are well-defined up to isomorphism. Every indecomposable direct summand of $M^\mu$ is isomorphic to some Young module and $Y^\mu\cong Y^\lambda$ if and only if $\mu=\lambda$. The Young modules have trivial sources. It is well known that if $\mu$ is a $p$-regular partition and $S^\mu$ is a simple Specht module then $S^\mu\cong Y^\mu$; see for example \cite[\S 1]{JCKMT}.

Every partition $\mu$ can be written as its $p$-adic expansion as follows. Let $\mu_{(0)}$ be the partition obtained from $\mu$ by successively stripping off all horizontal $p$-hooks. Then we have the coordinate-wise summation $\mu=p\mu(1)+\mu_{(0)}$ for some partition $\mu(1)$. Inductively, we can write $\mu(k-1)=p\mu(k)+\mu_{(k-1)}$ for $k\geq 2$ and hence $\mu=p^k\mu(k)+\cdots+p\mu_{(1)}+\mu_{(0)}$. The process ends if for some $k$, $\mu(k)$ has no removable horizontal $p$-hooks. Let $\mu_{(k)}=\mu(k)$. Then $\mu=p^k\mu_{(k)}+\cdots+p\mu_{(1)}+\mu_{(0)}$ is the $p$-adic expansion of $\mu$. Let $\rho(\mu)$ be the partition $$(\underbrace{p^k,\ldots,p^k}_{\text{$|\mu_{(k)}|$ factors}}\,,\ldots,\underbrace{p,\ldots,p}_{\text{$|\mu_{(1)}|$ factors}}\,,\underbrace{p^0,\ldots,p^0}_{\text{$|\mu_{(0)}|$ factors}})$$ of $n$ and $\mathfrak{S}_{\rho(\mu)}$ be the corresponding Young subgroup of $\sym{n}$.

\begin{thm}[{\cite{JGrab}, \cite[\S 4.1]{KE}}]\label{T: vtx Young} Any vertex of the Young module $Y^\mu$ is conjugate to a Sylow $p$-subgroup of $\mathfrak{S}_{\rho(\mu)}$.
\end{thm}

The main purpose of this article is to study the vertices and complexity of Specht modules. Thus, by Theorem \ref{P:complexity of mod} (i), it is crucial to understand the maximal elementary abelian $p$-subgroups of a symmetric group up to conjugation. Let $m$ be a positive integer. Let $(\mathbb{Z}_p)^m$ act on itself by the left regular action. This induces an injective group homomorphism $(\mathbb{Z}_p)^m\hookrightarrow \sym{p^m}$. Since the action is faithful we write $V_m(p)$ for the image of this homomorphism.

\begin{thm}[{\cite[\S VI Theorem 1.3]{AARM}}]\label{T:max ele abl} Let $n$ be a positive integer. There is a one-to-one correspondence between the ways of writing $n=i_0+i_1p+\cdots+i_rp^r$ for some non-negative integer $r$ such that all $i_0,\ldots,i_r$ are non-negative integers and $0\leq i_0\leq p-1$, and the set of representatives of all maximal elementary abelian $p$-subgroups of $\sym{n}$ up to conjugation, given by $$\underbrace{V_1(p)\times\cdots\times V_1(p)}_{\text{$i_1$ times}}\times \cdots\times \underbrace{V_r(p)\times\cdots\times V_r(p)}_{\text{$i_r$ times}}\hookrightarrow$$ $$\underbrace{\sym{p}\times\cdots\times \sym{p}}_{\text{$i_1$ times}}\times \cdots\times \underbrace{\sym{p^r}\times\cdots\times \sym{p^r}}_{\text{$i_r$ times}}\hookrightarrow \sym{n}.$$
\end{thm}

Note that the $p$-rank of the maximal elementary abelian $p$-subgroup of $\sym{n}$ given in Theorem \ref{T:max ele abl} is $i_1+2i_2+\cdots+ri_r$. We shall draw an easy conclusion, following Theorem \ref{P:Specht mod} (iv) and Theorem \ref{P:complexity of mod} (ii), that the complexity of a Specht module $S^\mu$ is bounded above by the $p$-weight of $\mu$. The statement remains valid if we replace the Specht module by the simple module $D^\mu$, whenever $\mu$ is $p$-regular, of $F\sym{n}$ as there is an analogous statement to Theorem \ref{P:Specht mod} (iv) for simple modules, though this is not of the main concern of this article.

\section{Main results}\label{S: main results}

We present our main results in this section.

\begin{thm}\label{T: main result 1} Let $F$ be a field of characteristic $p>0$ and $S^\mu$ be an indecomposable Specht module. Suppose that $S^\mu$ has an abelian vertex of $p$-rank $m$. Then $\mu$ is a $p^2$-core and the complexity of $S^\mu$ is $m$.
\end{thm}

\begin{thm}\label{T: main result 2} Let $F$ be a field of characteristic $p>0$ and $\mu$ be a partition. Suppose that $p\geq 3$, or $p=2$ and $\mu$ is $2$-regular. Suppose that the Specht module $S^\mu$ has an abelian vertex $Q$ of $p$-rank $m$. Let $c$ be the complexity of the Specht module $S^\mu$ and $w$ be the $p$-weight of the partition $\mu$. Then $c=m=w$ and $Q$ is conjugate to the elementary abelian $p$-subgroup $$\underbrace{V_1(p)\times \cdots\times V_1(p)}_{\text{$w$ factors}}\hookrightarrow \sym{p}\times \cdots \times \sym{p}\hookrightarrow \sym{n}.$$
\end{thm}

\begin{thm}\label{T: main result 3} Let $p=2$, $\mu$ be a $2$-regular partition and $Q$ be a vertex of the Specht module $S^\mu$. Let $w$ be a non-negative integer. Then the following statements are equivalent.
\begin{enumerate}
  \item [(i)] The vertex $Q$ is abelian of $2$-rank $w$.
  \item [(ii)] The vertex $Q$ is elementary abelian of $2$-rank $w$.
  \item [(iii)] The partition $\mu$ is a $4$-core of $2$-weight $w$.
\end{enumerate} In any of these cases, the Specht module $S^\mu$ has trivial source, complexity $w$ and is a simple Young module.
\end{thm}

\begin{rem} The hypothesis of Theorem \ref{T: main result 3} cannot be loosened in the following sense. For $p=2$, it is known that the Specht modules $S^{(n-r,1^r)}$ are indecomposable when $2\mid n$ \cite[Theorem 4.1]{GM}. In this case, the vertices of $S^{(n-r,1^r)}$ are the Sylow $2$-subgroups of $\sym{n}$ \cite[Theorem 4.5]{GMMP}. Thus, for instance, the Specht module $S^{(4,1,1)}$ has a non-abelian vertex but yet the partition $(4,1,1)$ is a $4$-core. So, in Theorem \ref{T: main result 3}, we cannot replace the $2$-regularity condition merely by the indecomposability condition.
\end{rem}

\section{Proof of the main results}\label{S: proofs}

In this section we shall prove our main theorems. We first prove Theorem \ref{T: main result 1}. We then divide the latter part of this section into two subsections. The first part concerns the proof of Theorem \ref{T: main result 2} when $p\geq 3$. The second part concerns the proofs of Theorem \ref{T: main result 2} when $p=2$ and Theorem \ref{T: main result 3}.

We state two known results upon which our proofs rely heavily.

\begin{thm}[{\cite{AP}, \cite[Theorem 1]{GK}}]\label{P:injection} Let $n_1,\ldots,n_s$ be positive integers. If there is an injective group homomorphism $\mathbb{Z}_{p^{n_1}}\times \ldots\times \mathbb{Z}_{p^{n_s}}\hookrightarrow \sym{m}$ then $p^{n_1}+\ldots+p^{n_s}\leq m$.
\end{thm}

For any positive integer $n$ and prime number $p$, we write $n_p$ for the prime power $p^a$ such that $p^a\mid n$ and $\gcd(p,n/p^a)=1$.

\begin{thm}[{\cite[Theorem 1.1]{CB}}]\label{T:Bessenrodt} Let $M$ be a $FG$-module for some finite group $G$. Suppose that $M$ is relatively $H$-projective for some subgroup $H$ of $G$. Let $Q\cong \mathbb{Z}_{p^{n_1}}\times \ldots\times \mathbb{Z}_{p^{n_m}}$ be an abelian subgroup of $H$ of $p$-rank $m$ with $n_1\geq \cdots \geq n_m\geq 1$. Let $c=c_Q(M)$ be the complexity of the $FQ$-module $\Res^G_Q M$. Then $$\left . \frac{|G:H|_p|Q|_p}{p^{n_1}\cdots p^{n_c}}\,\,\right |\,\, (\dim_F M)_p.$$
\end{thm}

We can now prove Theorem \ref{T: main result 1}.

\begin{proof}[Proof of Theorem \ref{T: main result 1}] We assume the notation as in Theorem \ref{T: main result 1}. Let $n=|\mu|$ and $w$ be the $p$-weight of $\mu$. Let $D$ be a defect group of the block $e_{\widetilde{\mu}}$ and $Q$ be a vertex of $S^\mu$. Without loss of generality, we may assume that $$Q\subseteq D\subseteq \sym{pw}.$$

We apply Theorem \ref{T:Bessenrodt} and Proposition \ref{L:complx rel vtx} by taking $M=S^\mu$, $H=Q\cong \mathbb{Z}_{p^{n_1}}\times \ldots\times \mathbb{Z}_{p^{n_m}}$ with $n_1\geq \cdots\geq n_m\geq 1$. Let $c=c_Q(S^\mu)=c_{\sym{n}}(S^\mu)$. We have $$p^{a-n_1-\cdots-n_m}\cdot p^{n_{c+1}}\cdots p^{n_m}=p^{a-n_1-\cdots-n_c}\mid (\dim_F S^\mu)_p$$ where $p^a=(n!)_p$. By Theorem \ref{P:Specht mod} (i) and the hook formula, we let $s$ be the non-negative integer such that $$p^s=\frac{p^a}{p^w\cdot (\dim_F S^\mu)_p}=\frac{\prod_{(i,j)\in [\mu]} (h_{i,j})_p}{p^w}.$$ Note that $s=0$ if and only if $\mu$ is a $p^2$-core; namely every hook length of $[\mu]$ is not divisible by $p^2$. With this notation, we deduce that
\begin{equation}\label{Eq:1}
s+w\leq n_1+\cdots+n_c.
\end{equation} On the other hand, we have the injection $Q\hookrightarrow\sym{pw}$. As such, by Theorem \ref{P:injection}, we get $p^{n_1}+\cdots+p^{n_m}\leq pw$; namely $p^{n_1-1}+\cdots+p^{n_m-1}\leq w$. Combining this inequality with (\ref{Eq:1}), we have \begin{equation}\label{Eq:2}
(p^{n_1-1}-n_1)+\cdots+(p^{n_c-1}-n_c)+p^{n_{c+1}-1}+\cdots+p^{n_m-1}+s\leq 0.
\end{equation} Since $p^{k-1}-k\geq 0$ for every positive integer $k$ we conclude that $m=c$, $s=0$, $p^{n_i-1}=n_i$ for all $1\leq i\leq c$ and $w=n_1+\cdots+n_c$. We have proved Theorem \ref{T: main result 1}.
\end{proof}

\subsection{The case where $p\geq 3$}\label{sS: p geq 3 case}

Assuming all the computations what we have just done in the proof of Theorem \ref{T: main result 1}, we continue to prove Theorem \ref{T: main result 2} for the case of $p\geq 3$.

\begin{proof}[Proof of Theorem \ref{T: main result 2} when $p\geq 3$] If $p$ is odd then, from $p^{n_i-1}=n_i$, necessarily $n_i=1$ for all $1\leq i\leq c$; namely $m=c=w$ and $Q$ is an elementary abelian $p$-group of $p$-rank $w$. Since $Q$ is also a subgroup of $\sym{pw}$, by Theorem \ref{T:max ele abl}, we let $$E=(V_1(p))^{i_1}\times \cdots\times (V_r(p))^{i_r}$$ be an maximal elementary $p$-subgroup of $\sym{pw}$ containing $Q$. Comparing their $p$-ranks, we have $w\leq i_1+2i_2+\cdots+ri_r$. On the other hand, we also have $i_1p+i_2p^2+\cdots+i_rp^r=pw$. Thus we deduce that $$i_2(p-2)+\cdots+i_r(p^{r-1}-r)\leq 0$$ and hence $i_2=\cdots=i_r=0$, $i_1=w$. Namely, $Q$ is necessarily conjugate to $(V_1(p))^w$. This completes the proof of Theorem \ref{T: main result 2} when $p$ is odd.
\end{proof}

\begin{rem} The author would like to point out that the proofs given above, for Theorem \ref{T: main result 1} and Theorem \ref{T: main result 2} when $p\geq 3$, is very similar to the proof of \cite[Theorem 1]{MW} in which Wildon dealt with the case where $m=1$.
\end{rem}

\subsection{The case where $p=2$}\label{sS: p=2 case}

Note that the proof of Theorem \ref{T: main result 1} assumes nothing about the $2$-regularity of the partition $\mu$, in the case of $p=2$. It works for all indecomposable Specht modules with abelian vertices. In view of Theorem \ref{T: main result 2}, we still have a result for $p=2$ case without assuming the $2$-regularity condition.

\begin{prop}\label{P: factors of vertices} Let $p=2$ and $S^\mu$ be indecomposable. Suppose that a vertex $Q$ of $S^\mu$ is abelian of $p$-rank $m$. Let $c$ be the complexity of $S^\mu$ and $w$ be the $2$-weight of $\mu$. Then $m=c$, $c\leq w\leq 2c$ and the direct factors of $Q$ are either $\mathbb{Z}_2$ or $\mathbb{Z}_4$.
\end{prop}
\begin{proof} By (\ref{Eq:2}) with $p=2$, we have $m=c$, $2^{n_i-1}=n_i$ for all $1\leq i\leq c$ and $w=n_1+\cdots+n_c$. Thus each of the $n_i$ is either $1$ or $2$.
\end{proof}

In general, we do not know when a Specht module over an even characteristic field is decomposable. Since vertices are defined only for indecomposable modules, in view of Theorem \ref{P:Specht mod} (iii), the $2$-regularity assumption is technically required.

Note that the proof which we have just given in \S \ref{sS: p geq 3 case} for Theorem \ref{T: main result 2} fails when $p=2$ because of the possibility that $Q$ may have $\mathbb{Z}_4$ as a direct factor. Thus we need to take a closer look at these partitions $\mu$ which are $2$-regular and the vertices of $S^\mu$ are abelian. 

\begin{lem}\label{L: 4-core} Let $\mu=(\mu_1,\ldots,\mu_r)$ be a partition. Then $\mu$ is both a $2$-regular and $4$-core partition if and only if exactly one of the following holds.
\begin{enumerate}
  \item [(i)] $\mu_r=1$ and there exists a number $0\leq s\leq r-1$ such that $\mu_i-\mu_{i+1}=3$ for all $1\leq i\leq s$ and $\mu_i-\mu_{i+1}=1$ for all $s+1\leq i\leq r-1$.

  \item [(ii)] $\mu_r\in \{2,3\}$ and $\mu_i-\mu_{i+1}=3$ for all $1\leq i\leq r-1$.
\end{enumerate}
\end{lem}
\begin{proof} Suppose that $\mu$ is $2$-regular and is a $4$-core. It is easy to see that $1\leq \mu_i-\mu_{i+1}\leq 3$ for all $1\leq i\leq r-1$ and $\mu_i-\mu_{i+1}\neq 2$ for all $1\leq i\leq r-2$. Suppose that $\mu_j-\mu_{j+1}=3$ and $\mu_{j-1}-\mu_j=1$ for some $2\leq j\leq r-1$. Then the node $(j-1,\mu_{j-1}-2)$ has hook length $4$. If $\mu_r$ is $2$ or $3$, then $\mu_{r-1}-\mu_r\neq 1$; otherwise the node $(r-1,1)$, respectively $(r-1,2)$, has hook length $4$. This proves one direction of the characterization. The converse is easy.
\end{proof}

%Let us consider some small examples. Let $\mu$ be either $(5,2)$, $(6,3)$, $(8,5,2)$ or $(7,4,3,2,1)$. Note that these Specht modules are simple, for example see \cite[Appendix]{GJ1} and \cite{JM}. In these examples, $S^\mu\cong D^\mu$. The simple modules, for the first two examples, has vertices conjugate to $V_1(p)\times V_1(p)\times V_1(p)$ as computed in \cite[Table A.1]{SDBKRZ} as a subgroup of $\sym{7}$ and $\sym{9}$ respectively. The examples provide evidence to convince the author that Theorem \ref{T: main result 2} should work for the case $p=2$. Indeed, the Specht modules of the partitions which we have characterized in Lemma \ref{L: 4-core} are simple. The characterization of simple Specht modules when $p=2$ has been done in \cite[Theorem]{GJAM}.

Using the above characterization of $2$-regular and $4$-core partitions, and the characterization of simple Specht modules for $p=2$ \cite[Main Theorem]{GJAM}, we obtain the following corollary.

\begin{cor}[{of \cite[Main Theorem]{GJAM}}]\label{C: 4-core is simple} Let $p=2$ and $\mu$ be a $2$-regular partition. If $\mu$ is a $4$-core then $S^\mu$ is simple and $D^\mu\cong S^\mu\cong Y^\mu$.
\end{cor}
%\begin{proof} We follow the notation of \cite{GJAM}. Let $l(k)$ be the least non-negative integer such that $k<2^{l(k)}$, for any integer $k$. We shall use \cite[Theorem (i)]{GJAM}; namely if $\mu_i-\mu_{i+1}\equiv -1(\mod 2^{l(\mu_{i+1}-\mu_{i+2})})$ for all $i\geq 1$ then $S^\mu$ is simple. With respect to Lemma \ref{L: 4-core}, we consider two cases.
%\begin{enumerate}
%  \item [(i)] In this case $\mu_r=1$. For $i\leq s-1$, $\mu_i-\mu_{i+1}=3\equiv -1(\mod 2^{l(\mu_{i+1}-\mu_{i+2})=2})$. For $s\leq i\leq r-1$, $\mu_s-\mu_{s+1}=3\equiv -1(\mod 2^{l(\mu_{i+1}-\mu_{i+2})=1})$. For $i=r$, $\mu_r-\mu_{r+1}=1\equiv -1(\mod 2^{l(\mu_{r+1}-\mu_{r+2})=0})$.

%  \item [(ii)] In this case $\mu_r=2,3$. For $i\leq r-2$, $\mu_i-\mu_{i+1}=3\equiv -1(\mod 2^2)$. For $i=r-1$, for both cases, $\mu_{r-1}-\mu_r=3\equiv -1(\mod 2^2)$. For $i=r$, $\mu_{r}-\mu_{r+1}\equiv -1(\mod 2^0)$.
%\end{enumerate} Simple Specht module $S^\mu$ with $\mu$ $p$-regular is isomorphic to the Young module $Y^\mu$.
%\end{proof}

\begin{cor}\label{C: vtx 4-core} Let $p=2$ and $\mu$ be a $2$-regular partition of $n$. Suppose that $\mu$ is a $4$-core. Then any vertex $Q$ of $S^\mu$ is conjugate to $$\underbrace{V_1(2)\times \cdots\times V_1(2)}_{\text{$w$ factors}}\hookrightarrow (\sym{2})^w\hookrightarrow \sym{n}$$ where $w$ is the $2$-weight of $\mu$. Furthermore, the complexity of $S^\mu$ is $w$.
\end{cor}
\begin{proof} Let $\mu=2^k\mu_{(k)}+\cdots+2\mu_{(1)}+\mu_{(0)}$ be the $2$-adic expansion of $\mu$ with $|\mu_{(k)}|\neq 0$. We claim that $0\leq k\leq 1$ and $\mu_{(0)}$ is the $2$-core of $\mu$. Note that there are only two types of removable skew $2$-hooks, the horizontal $2$-hook and the vertical $2$-hook. Since $\mu_{(0)}$ is obtained from $\mu$ by successively removing all horizontal $2$-hooks it is then not a $2$-core if there is some removable vertical $2$-hook of $\mu_{(0)}$. In this case, the original partition $\mu$ has two successive non-empty rows whose difference of their sizes is divisible by $2$. However, this is not allowed by the characterization in Lemma \ref{L: 4-core}. Thus $\mu_{(0)}=\widetilde{\mu}$. It is clear that, in general, if $\mu$ is a $p^k$-core for some positive integer $k$ then $\mu_{(l)}=\varnothing$ for all $l\geq k$.

By Corollary \ref{C: 4-core is simple}, the Specht module $S^\mu$ is isomorphic to the Young module $Y^\mu$. By Theorem \ref{T: vtx Young}, we deduce that any vertex $Q$ of the Specht module $S^\mu\cong Y^\mu$ is conjugate to a Sylow $2$-subgroup of $\mathfrak{S}_{\rho(\mu)}=(\sym{2})^{|\mu_{(1)}|}$ where $$|\mu_{(1)}|=(n-|\mu_{(0)}|)/2=(n-|\widetilde{\mu}|)/2=w.$$ Since the Young module $Y^\mu$ has a trivial source; namely the $FQ$-module $F$, by Proposition \ref{L:complx rel vtx} (ii), we have $$c_{\sym{n}}(S^\mu)=c_Q(F)=w.$$ The proof of our corollary is now complete.
\end{proof}

The proofs for Theorem \ref{T: main result 2} when $p=2$ and Theorem \ref{T: main result 3} are now clear:

\begin{proof}[Proof of Theorem \ref{T: main result 2} when $p=2$] We assume all the notation in Theorem \ref{T: main result 2}. By Theorem \ref{T: main result 1}, the partition $\mu$ is necessarily a $4$-core. By Corollary \ref{C: vtx 4-core}, we get $m=w=c$ and $Q$ is conjugate to $(V_1(2))^w$.
\end{proof}

\begin{proof}[Proof of Theorem \ref{T: main result 3}] Corollary \ref{C: vtx 4-core} shows that (iii) implies (ii). The implication from (ii) to (i) is trivial. Theorem \ref{T: main result 1} shows that (i) implies (iii). The final statement comes from Corollaries \ref{C: 4-core is simple} and \ref{C: vtx 4-core}.
\end{proof}

\section{Some applications of the main results}\label{S: consequences}

In \cite[Theorem 1]{MW}, M. Wildon showed that the vertices of an indecomposable Specht module $S^\mu$ are non-trivial cyclic if and only if the $p$-weight of $\mu$ is $1$. We have a generalization of this result.

\begin{cor}\label{C: abl vtx in abl dft} Let $p$ be the characteristic of the field $F$ and suppose that $p\geq 3$. Let $1\leq m\leq p-1$ be an integer and $\mu$ be a partition such that the Specht module $S^\mu$ is indecomposable. Then a vertex $Q$ of $S^\mu$ is abelian of $p$-rank $m$ if and only if the $p$-weight of $\mu$ is $m$. In this case, $Q$ is necessarily elementary abelian.
\end{cor}
\begin{proof} Assume that $Q$ is abelian. Theorem \ref{T: main result 2} shows that $Q$ is elementary abelian and $m=w$ where $w$ is the $p$-weight of $\mu$. Conversely, suppose that $m=w$. In this case, we have the abelian defect case. Thus $Q$ is necessarily abelian as a subgroup of some defect group of $S^\mu$. By Theorem \ref{T: main result 2}, $Q$ is has $p$-rank $w=m$.
\end{proof}

\begin{rem} By virtue of \cite[Theorem 1]{MW}, Corollary \ref{C: abl vtx in abl dft} also holds for $p=2$.
\end{rem}

There is a special class of Specht modules $S^\mu$ whose partition $\mu$ is of the form $({\mu_1}^{a_1},{\mu_2}^{a_2},\ldots,{\mu_r}^{a_r})$ where both $\mu_i,a_i$ are multiples of $p$ for every $1\leq i\leq r$. We call them the $(p\times p)$-partitions. As another application of our results, we study the vertices of the Specht modules $S^\mu$ where $\mu$ is a $(p\times p)$-partition for $p\geq 3$. We have seen in Lemma \ref{L:vtx fld ext} that the vertices of an indecomposable module remain unchanged, should the module remains indecomposable upon field extension. Since all Specht modules in the case of $p\geq 3$ are indecomposable, we may assume that $F$ is algebraically closed for the rest of this section.

Let $E_p=V_1(p)\times \cdots\times V_1(p)$ ($p$ factors) be the maximal elementary abelian $p$-subgroup of $\sym{p^2}$ as in Theorem \ref{T:max ele abl}. We have $c_{\sym{p^2}}(S^{(p^p)})=c_{E_p}(S^{(p^p)})=p-1$ \cite[Theorem 3.1(i)]{KJL}. Furthermore, for $p\geq 3$, with respect to the generators $$(1,\ldots,p),(p+1,\ldots,2p),\ldots,(p^2-p+1,\ldots,p^2)$$ of $E_p$, a description of the union of irreducible components $W_{p-1}$ of dimension $p-1$ of $V^\#_{E_p}(S^{(p^p)})$ and the vanishing ideal of $W_{p-1}$ is given in \cite[Theorem 3.1(ii)]{KJL}. The vanishing ideal $I$ of $W_{p-1}$ is generated by the element $$f(x_1,\ldots,x_p)=(x_1\cdots x_p)^{p-1}\widetilde{f}+\sum_{i=1}^p {x_1}^{n(p-1)}\cdots \widehat{{x_i}^{n(p-1)}}\cdots {x_p}^{n(p-1)}$$ for some homogeneous polynomial $\widetilde{f}$ and positive integer $n$, and where the term ${x_1}^{n(p-1)}\cdots \widehat{{x_i}^{n(p-1)}}\cdots {x_p}^{n(p-1)}$ is the product of all ${x_j}^{n(p-1)}$'s with $1\leq j\leq p$ except the term ${x_i}^{n(p-1)}$. In particular, we have the following lemma.

\begin{lem}\label{L: irreducible component of (p^p)} Any irreducible component of $W_{p-1}\subseteq V^\#_{E_p}(S^{(p^p)})$ is not a hyperplane of $\mathbb{A}^p(F)$.
\end{lem}
\begin{proof} Suppose that $W_{p-1}$ contains a hyperplane, i.e. $f(x_1,\ldots,x_p)\in I(W_{p-1})\subseteq \langle a_1x_1+\cdots+a_px_p\rangle$ for some non-zero vector $v=(a_1,\ldots,a_p)\in \mathbb{A}^p(F)$. Let $g\in F[x_1,\ldots,x_p]$ be the polynomial such that
\begin{equation}\label{Eq:3}
f(x_1,\ldots,x_p)=(a_1x_1+\cdots+a_px_p)g(x_1,\ldots,x_p).
\end{equation} Suppose that at least two of the coordinates of $v$ are non-zero; say $a_r$ and $a_s$. Let $1\leq i\leq p$ be such that $i$ is different from either $r$ or $s$. Then there exists elements $b_1,\ldots,b_p$ of $F$ with $b_i=0$, and $b_j\neq 0$ for all $1\leq j\leq p$ and $j\neq i$, such that $a_1b_1+\cdots+a_pb_p=0$. Substitute the values $b_1,\ldots,b_p$ into Equation (\ref{Eq:3}), we have the contradiction where $$0\neq {b_1}^{n(p-1)}\cdots \widehat{{b_i}^{n(p-1)}}\cdots {b_p}^{n(p-1)}=0\cdot g(b_1,\ldots,b_p)=0.$$ This shows that the coordinates of $v$ are all zero except $a_j\neq 0$ for some unique $1\leq j\leq p$. We deduce that $f(x_1,\ldots,x_p)=a_jx_jg(x_1,\ldots,x_p)$ and hence $x_j$ is a factor of $f(x_1,\ldots,x_p)$. However, it is clear from the description of $f$ that $x_j$ is not a factor. The proof is now complete.
\end{proof}

For $p\geq 3$, Hemmer show that the complexity of the Specht module $S^\mu$ when $\mu$ is a $(p\times p)$-partition is strictly less than the $p$-weight of $\mu$ \cite[Corollary 1.4]{DH}. Thus the following corollary is immediate by applying Theorem \ref{T: main result 2}.

\begin{cor}\label{C: p by p vtx is non-abelian} Suppose that $p\geq 3$ and $\mu$ is a $(p\times p)$-partition. Then the vertices of $S^\mu$ are non-abelian.
\end{cor}

The next corollary describes the vertices of the Specht module $S^{(p^p)}$ for $p\geq 3$.

\begin{cor}\label{C: p by p vtx} For $p\geq 3$, any vertex of $S^{(p^p)}$ is a Sylow $p$-subgroup of $\sym{p^2}$.
\end{cor}
\begin{proof} Let $Q$ be a vertex of $M=S^{(p^p)}$ and $S$ be a source of $M$. By Mackey Decomposition Theorem, we have $$\Res^{\sym{p^2}}_{E_p}M\mid \Res^{\sym{p^2}}_{E_p}\Ind^{\sym{p^2}}_{Q}S\cong \bigoplus_{g\in E_p\backslash \sym{p^2}/Q} \Ind^{E_p}_{{}^g Q\cap E_p} \Res^{{}^g Q}_{{}^gQ\cap E_p}{}^g S$$ where $E_p\backslash \sym{p^2}/Q$ is a set of double coset representatives of $(E_p,Q)$ in $\sym{p^2}$. Suppose that ${}^gQ\cap E_p\neq E_p$ for all $g$; namely ${}^gQ\cap E_p$ is a proper subgroup of $E_p$ and hence has order at most $p^{p-1}$. Let $$U_g=V^\#_{E_p}(\Ind^{E_p}_{{}^g Q\cap E_p} \Res^{{}^g Q}_{{}^gQ\cap E_p}{}^g S).$$ In fact, $\dim U_g=\dim V^\#_{{}^gQ\cap E_p} (\Res^{{}^g Q}_{{}^gQ\cap E_p}{}^g S)$. Since $U_g$ is an induced module induced from a proper subgroup of $E_p$ we have that the rank variety $$U:=V^\#_{E_p}(\Res^{\sym{p^2}}_{E_p}\Ind^{\sym{p^2}}_{Q}S)=\bigcup U_g\subseteq \mathbb{A}^p(F)$$ is a finite union of subvarieties of dimension at most $p-1$. If $U_g$ has dimension $p-1$ for some $g$ then it is necessarily that $V^\#_{{}^gQ\cap E_p} (\Res^{{}^g Q}_{{}^gQ\cap E_p}{}^g S)\cong \mathbb{A}^{p-1}(F)$ and hence $U_g$ is a union of hyperplanes of $\mathbb{A}^p(F)$.

Let $N$ be a $kE_p$-module such that $N\oplus \Res^{\sym{p^2}}_{E_p}M=\Res^{\sym{p^2}}_{E_p}\Ind^{\sym{p^2}}_{Q}S$. We have $V^\#_{E_p}(M)\cup V^\#_{E_p}(N)=U$. Let $V$ be an irreducible component of $V^\#_{E_p}(M)$ of dimension $p-1$. By the unique decomposition property of varieties into their irreducible varieties, see \cite[Corollary 1.6]{RH}, we conclude that $V$ is a hyperplane of $\mathbb{A}^p(F)$. This contradicts to Lemma \ref{L: irreducible component of (p^p)} and hence we conclude that ${}^gQ\cap E_p=E_p$ for some $g\in \mathfrak{S}_{p^2}$; namely $E_p\subseteq Q$ up to some conjugation.

Since $S^{(p^p)}$ has complexity $p-1$, different from the $p$-weight of $(p^p)$, by Theorem \ref{T: main result 2} or Corollary \ref{C: p by p vtx is non-abelian}, the vertices cannot be abelian. Thus $Q$ has order strictly larger than $p^p$, i.e. $p^{p+1}$ and hence is a Sylow $p$-subgroup of $\sym{p^2}$.
\end{proof}

\begin{rem} The case of $p=2$ is dealt with in \cite[Lemma 6]{MW}: the Specht module $S^{(2,2)}$ has the Klein four group $V_2(2)$ as its vertex, which is different from the description in Corollary \ref{C: p by p vtx}.
\end{rem}

\section{Some further questions}\label{S: some further questions}

In this section, we include some further questions. The following question is natural following our Theorem \ref{T: main result 3}.

\begin{ques}\label{Q: Question 1}
Given that $\mu$ is $p$-regular and $p\geq 3$, are the vertices of the Specht module $S^\mu$ abelian if $\mu$ is a $p^2$-core?
\end{ques}

%Corollary \ref{C: p by p vtx is non-abelian} and \ref{C: p by p vtx} provide some examples of Specht modules $S^\mu$ where the partition $\mu$ is a $p^2$-core but yet the vertices are non-abelian.

To answer Question \ref{Q: Question 1}, our method used in \S \ref{sS: p=2 case} does not help so much. The main obstruction is that there are Specht modules $S^\mu$ which are not Young modules, with $\mu$ $p$-regular and a $p^2$-core (for example, the Specht module $S^{(6,3)}$ with $p=3$). Readers may have observed that, in order to achieve Theorem \ref{T: main result 3}, the main idea is to prove the simplicity of the corresponding Specht modules as in Corollary \ref{C: 4-core is simple}. Thus we can partly answer Question \ref{Q: Question 1} by imposing this crude assumption.

\begin{prop}\label{P: classfication for simple, p odd} Let $p\geq 3$, $\mu$ be a $p$-regular partition and $Q$ be a vertex of the Specht module $S^\mu$. Suppose that $S^\mu$ is simple. Let $w$ be a non-negative integer. Then the following statements are equivalent.
\begin{enumerate}
  \item [(i)] The vertex $Q$ is abelian of $p$-rank $w$.
  \item [(ii)] The vertex $Q$ is elementary abelian of $p$-rank $w$.
  \item [(iii)] The partition $\mu$ is a $p^2$-core of $p$-weight $w$.
\end{enumerate} In any of these cases, the Specht module $S^\mu$ has trivial source, complexity $w$ and is a simple Young module.
\end{prop}
\begin{proof} We only need to show that (iii) implies (i). As we have mentioned in the proof of Corollary \ref{C: vtx 4-core}, the $p$-adic expansion of $\mu$ is $\mu_{(0)}+p\mu_{(1)}$, here we allow $\mu_{(1)}=\varnothing$. By the simplicity condition, we deduce that $S^\mu\cong Y^\mu$ and hence $Q$ is conjugate to a Sylow $p$-subgroup of $\sym{\rho(\mu)}=(\sym{p})^{|\mu_{(1)}|}$ by Theorem \ref{T: vtx Young}. So $Q$ is abelian, necessarily of $p$-rank $w$ and the Specht module $S^\mu$ has complexity $w$ by Theorem \ref{T: main result 1}.
\end{proof}

Of course, ultimately, we will be delighted to achieve the following goal, if the computation of the vertices of all Specht modules is difficult.

\begin{ques}
Classify all indecomposable Specht modules with abelian vertices.
\end{ques}

As we have pointed out in the beginning of \S \ref{sS: p=2 case}, the technical assumption of $2$-regularity is required in Theorem \ref{T: main result 2}. The result will be nicer if it works for all indecomposable Specht modules with abelian vertices. The possibility of the existence of an example of a Specht module $S^\mu$ with abelian vertices $Q$ which contain  $\mathbb{Z}_4$ as a direct factor is the main obstruction.We do not have an example of a Specht module that support this possibility of Proposition \ref{P: factors of vertices}.

\begin{ques} For $p=2$, is there an indecomposable Specht module $S^\mu$ whose vertices are abelian and contain $\mathbb{Z}_4$ as a direct factor?
\end{ques}

\subsection*{Acknowledgment} The results of this article have been partly achieved during my days in Aberdeen \cite[Remark 4.14]{KJL2} and partly done in Singapore. I would like to take the chance to thank Dave Benson and Radha Kessar. I would like to thank Kai Meng Tan for his availability for discussion and careful reading of this article, and Mark Wildon for pointing out Theorem \ref{P:injection} and improving an earlier version of Theorem \ref{T: main result 1}.

I would like to thank the referees for their valuable suggestions.

\end{document}